\documentclass[10pt, reqno]{amsart}
\usepackage{amssymb,latexsym}
\newtheorem*{proposition*}{}
\newtheorem*{LemA}{Lemma A}
\newtheorem*{LemB}{Lemma B}
\newtheorem*{LemC}{Lemma C}
\newtheorem*{Lem10.4}{Lemma 10.4}
\newtheorem*{Lem17.3}{Lemma 17.3}
\newtheorem*{Lem18.3}{Lemma 18.3}
\newtheorem*{Lem18.5}{Lemma 18.5}
\newtheorem*{Lem19.0}{Lemma 19.0}
\newtheorem*{Lem19.1}{Lemma 19.1}
\newtheorem*{Lem19.1.3}{Lemma 19.1.3}

\newtheorem*{ThA}{Theorem A}
\newtheorem*{ThB}{Theorem B}

\newcommand{\e}{{\varepsilon}}

\newcommand{\G}{{\Gamma}}
\newcommand{\D}{\Delta}

\newcommand{\p}{\partial}
\newcommand{\E}{{\textup{E}}}

\begin{document}

\title[]{Weakly finitely presented infinite periodic groups  }
\author{S. V. Ivanov}
\address{Department of Mathematics\\
University of Illinois \\
Urbana, \  IL 61801}
\email{ivanov@math.uiuc.edu} 
\thanks{Supported in part by NSF grants DMS 98-01500,  DMS 00-99612}
\begin{abstract}
A group $G$ given by a presentation 
$G = \langle \mathcal A \; \| \; \mathcal R \rangle$ is called
weakly finitely presented if every finitely generated  subgroup of $G$, 
generated by 
(images of) some words in $\mathcal A^{\pm 1}$, 
is naturally isomorphic to the subgroup of 
a group $G_0 = \langle \mathcal A_0 \; \| \; \mathcal R_0 \rangle$, 
where $\mathcal A_0 \subseteq \mathcal A$,   
$\mathcal R_0 \subseteq \mathcal R$ are finite,  
generated by  (images of) the same words.  
In the article, weakly finitely presented  periodic groups 
which are not locally finite are constructed. 
\end{abstract}

\maketitle


\section{Introduction}

Let a group $G$ be given by a presentation
\begin{equation}\label{1}
G = \langle  \mathcal A  \  \|\  \  R =1, \ R  \in \mathcal R  \rangle  , 
\end{equation}
where $\mathcal A$ is an alphabet, $\mathcal R$ is a set of defining relators (which are reduced words in the alphabet 
$\mathcal A^{\pm 1} =\mathcal A \cup \mathcal A^{-1}$).  
Let $F(\mathcal A)$ be a free group in  the alphabet $\mathcal A$, 
$N(\mathcal R) $ the normal closure of $\mathcal R$ in $F(\mathcal A)$ and 
$$
\psi_0 : F(\mathcal A) \to G = F(\mathcal A)/ N(\mathcal R)
$$
be the natural  homomorphism. 

This group presentation (1) is referred to as {\em weakly finite} if for every finite set 
$\mathcal S = \{ W_1, \dots, W_m\}$ of words $W_1, \dots, W_m$  in $\mathcal A^{\pm 1}$,
called $\mathcal A$-{\em words}, there are finite subsets 
$\mathcal A_{\mathcal S} \subseteq  \mathcal A$ and $\mathcal R_{\mathcal S} \subseteq  \mathcal R$ such that 
$W_1, \dots, W_m$  are $\mathcal A_{\mathcal S}$-words and the subgroup 
$\langle \mathcal S^{\psi_0} \rangle$ of $G = F(\mathcal A)^{\psi_0} $ generated by $\mathcal S^{\psi_0}$ 
is naturally 
isomorphic to  the subgroup $\langle \mathcal S^{\psi_{0, \mathcal S}} \rangle$ of 
the group $G_{\mathcal S} = F(\mathcal A_{\mathcal S})^{\psi_{0, \mathcal S}}$, where 
$$
\psi_{0, \mathcal S} :   F(\mathcal A_{\mathcal S}) \to G_{\mathcal S} = 
F(\mathcal A_{\mathcal S})/ N(\mathcal R_{\mathcal S})
$$
is the natural homomorphism. Note that, in the foregoing notation,  
there is always a  natural homomorphism 
$$
\langle \mathcal S^{\psi_{\mathcal S}} \rangle \to \langle \mathcal S^\psi \rangle$$ 
for the natural homomorphism $\langle \mathcal S \rangle \to \langle \mathcal S^\psi \rangle$ 
factors out through 
$\langle \mathcal S^{\psi_{\mathcal S}} \rangle$.

Accordingly, a group $G$ given by (1) is called  {\em weakly finitely presented} if (1) is
a weakly finite  presentation. 

Recall that  the presentation (1) is called {\em finite} if both $\mathcal A$ and $\mathcal R$ are finite.  It is easy to see that  if $\mathcal A$ is finite then (1) is  weakly finite if and only if there is a finite subset 
$\mathcal R_0 \subseteq \mathcal R$ such that  
$N(\mathcal R_0) = N(\mathcal R)$.  It  is also easy to see that 
the presentation (1) is weakly finite if every finitely generated subgroup of  $G$ is finitely presentable (that is, has a finite presentation). In particular, 
if $G$  is a locally finite group then $G$ is weakly finitely presented. 

It seems rather natural to ask whether every periodic  weakly finitely presented  
group is locally finite. This problem can be regarded  as a weakened version of a 
long-standing problem, attributed to P.S. Novikov, on the existence of  a finitely presented infinite periodic group. In this article, we will solve the former problem in the negative thus making some progress towards the  Novikov problem.  However, our construction
of weakly finitely presented  non-locally finite periodic groups does not even allow to bound the orders of elements  and the problem on the existence of  
weakly finitely presented   groups of bounded exponent that are not 
locally finite  might be worth further investigation. 

Recall that first examples of finitely generated infinite periodic groups of unbounded exponent  were constructed by Golod \cite{eG64}.  
Later more examples of such groups (which are rather easy to construct)
were  found by Aleshin \cite{svA72},  Grigorchuk \cite{rG80}, Gupta and Sidki \cite{GS83}, \cite{nG89}. 
Also, recall that the Burnside problem \cite{B02}  on periodic groups 
asks about the existence of   finitely generated  groups of {\em exponent} $n$ (that is, satisfying the identity $x^{n} \equiv 1$).  
For odd exponents $n \gg 1$ the Burnside  problem was solved 
by Novikov and Adian  \cite{NA68}  in 1968
(a simpler geometric solution for odd  $n \gg 1$  was later found by
Ol'shanskii's \cite{O82}, \cite{O89}). The case of even exponents in the Burnside problem  turned out to be much more difficult and required creation of new heavy machinery, see author's article \cite{I94} 
(or Lysenok's article \cite{L96}),  
that was aimed to handle  noncyclic centralizers (more  applications of this 
machinery  can be found in \cite{IO96}, \cite{IO97},  \cite{I00}). 
Curiously, many parts of the machinery of article \cite{I94}
will be "recycled"  in this paper in order to construct weakly finitely presented   
non-locally finite periodic groups  whose inductive construction also 
generates noncyclic centralizers. 

Recall that an $m$-generator free Burnside group $B(m,n)$  
of exponent $n$ is the quotient  $F_m/F_m^n$, where $F_m$ 
is a free group of rank $m$ and  $F_m^n$ is the subgroup of  
$F_m$ generated by all  $n$th powers of elements of $F_m$. 
By Novikov-Adian theorem (see \cite{NA68}, \cite{siA75}),
$B(m, n)$ is infinite if $m >1$ and $n$ is odd, $n \ge 665$ 
(see also \cite{O82}, \cite{O89}, \cite{I94}). 

\begin{ThA}\label{theor1} 
Suppose that $N_0 > 2^{16}$ is an odd integer. Then there is a 
weakly finitely presented  periodic group 
$$
G =  \langle \; \mathcal A \; \|  \; R=1, \; R \in \mathcal R   \; \rangle ,
$$ 
where both $\mathcal A$ and  $\mathcal R$ are countably infinite, 
such that $G$ has unbounded exponent and 
$G$ is not a locally finite group for $G$ contains
a subgroup isomorphic to a  
$2$-generator free Burnside group $B(2, N_0)$ of exponent $N_0$. 
\end{ThA}

A detailed description of our construction of weakly finitely presented  periodic groups 
is given in the following.

\begin{ThB}\label{theor2}
 Suppose that $N_0 > 2^{16}$ is an odd integer, 
$$
G_0 = \langle \; \mathcal A_0 \; \|  \; R=1, \; R \in \mathcal R_0 \; \rangle
$$
is a finitely presented group and the alphabet 
$\mathcal A_0$ contains at least two letters. 
Then one can construct  finite alphabets
$$
\mathcal A_0 \subset \mathcal A_1 \subset \mathcal A_2 \subset   \dots 
$$ 
and finite sets 
$$
\mathcal R_0 \subset  \mathcal R_1 \subset \mathcal R_2 \subset  \dots 
$$
of words in  $\mathcal  A_0^{\pm 1}$, $\mathcal  A_1^{\pm 1}$, 
$\mathcal  A_2^{\pm 1}, \dots$, respectively,  such that the following hold. 

Let 
$$
G_j = \langle \; \mathcal A_j \; \|  \; R=1, \ R \in \mathcal R_j \; \rangle ,
$$
$F(\mathcal A_j)$ be a free group in the alphabet  $\mathcal A_j$ and 
$$
\psi_{0,j} :  F(\mathcal A_j) \to G_j
$$ 
be the natural homomorphism, $j=0, 1,2,\dots$.
Then for all $j \geq 0$,  $k \ge 1$  it is true that
$$
F(\mathcal A_j)^{\psi_{0, j+k}} =  F(\mathcal A_j)^{\psi_{0, j}}/ 
( F(\mathcal A_j)^{\psi_{0, j}})^{N_0^{3^j}} = G_j / G_j^{N_0^{3^j}}  .
$$ 
In particular,  the group
$$
G_\infty = \langle \;  \mathcal A_\infty = \bigcup_{j=0}^\infty \mathcal A_j 
\; \|  \; R=1, \  R \in \mathcal R_\infty = \bigcup_{j=0}^\infty \mathcal R_j   \; \rangle 
$$
contains isomorphic copies of quotients 
$$
G_0/ G_0^{N_0}, \ \  G_1/ G_1^{N_0^3}, \ \  G_2 / G_2^{N_0^{3^2}}, \dots,  \ \ 
G_j / G_j^{N_0^{3^j}}, \dots
$$
which, except possibly for  $G_0/ G_0^{N_0^3}$, are not locally finite and if 
$W_1, \dots, W_m$ are words in 
$\mathcal A_j^{\pm 1}$ and 
$$
\psi_{0, \infty} :  F(\mathcal A_\infty) \to G_\infty
$$ 
is the natural homomorphism then the subgroup 
$\langle W_1^{\psi_{0, \infty}}, \dots, W_m^{\psi_{0, \infty}} \rangle$
of $G_\infty $  has exponent $N_0^{3^{j+1}}$ 
(i.e., it satisfies the identity $x^{N_0^{3^{j+1}}} \equiv 1$) 
and naturally isomorphic  to  
$$
\langle W_1^{\psi_{0, j+1}}, \dots, W_m^{\psi_{0, j+1}} \rangle \subseteq G_{j+1} .
$$
\end{ThB}

\section{The Main Construction}

Let  $ \mathcal A = \{ a_i \; | \;  i \in I \}$, 
where $I$ is a (finite) index set with $|I| >1$,  be an alphabet.
We also consider alphabets 
$$
\mathcal B  = \{ b_i \;  | \;  i \in I \} ,  \  \   
\mathcal X = \{ x_i \;  |  \; i \in I \} , \   \mbox{and}  \   \{ c,  y \}  
$$
 such  that  $\mathcal {A, B, X}, \{c , \; y \}$ are pairwise disjoint.

Suppose that $\mathcal R$ is a set of words in 
$\mathcal {A}^{\pm 1} = \mathcal {A} \cup \mathcal {A}^{-1}$ and 

\begin{equation}\label{2}
H = \langle  \mathcal {A}  \  \|\  \  R =1, \ R  \in \mathcal R  \rangle  
\end{equation}
is a group presentation whose set of generators is $\mathcal A$ and whose set of 
defining relators is $\mathcal R$. Consider the following defining relations

\begin{gather}
x_i c x_i^{-1} =  c b_i , \  \ i \in I ,   \label{3} \\ 
y b_i  y^{-1} =   b_i  a_i , \  \ i \in I ,    \label{4} \\ 
c^n =1 ,  \label{5} \\ 
x_i b_j  =  b_j  x_i  , \  \ i, j  \in I ,  \label{6} \\ 
a_i b_j  =   b_j a_i  , \  \ i, j  \in I , \label{7} \\ 
a_i  c  =   c a_i  , \  \ i   \in I ,  \label{8} \\ 
y c  =   c y , \  \ i   \in I .  \label{9}
\end{gather}

Let a group $G$ be given by  a presentation whose alphabet is 
$$
\mathcal U = \mathcal A \cup \mathcal B \cup \mathcal X \cup \{c, y\}
$$
and whose set of   
defining relations consists of  defining relations of (2) and relations (3)--(9), 
thus 
\begin{equation} \label{10}
G = \langle  \ a_i, b_i, x_i, c, y,    i \in I  \  \|\  \  R =1, \ R  \in \mathcal R, \  
(3)-(9)  \  \rangle .  
\end{equation}

Denote a free group in the alphabet $\mathcal A$ by $F(\mathcal A)$ and let 
$$
{\alpha_0} : F(\mathcal A) \to H ,  \quad   \beta_0 :  F(\mathcal U) \to G 
$$
be natural  homomorphisms. 

A word $W$ in the alphabet $\mathcal A^{\pm 1}$  
is called an $\mathcal A$-{\em word} and 
denoted  by $W = W(\mathcal A)$. A substitution $b_i \to a_i, i \in I$, 
turns an $\mathcal A$-word $W(\mathcal A)$ into a  $\mathcal B$-word denoted by $W(\mathcal B)$.

\begin{LemA} 
The map $a_i^{\alpha_0}  \to a_i^{\beta_0}$, $i \in I$, extends to a homomorphism 
$$
\psi_1 : H \to G
$$ 
whose kernel $\textup{Ker}\psi_1$ is 
$H^n = \langle h^n \; | \; h \in H \rangle$. 
\end{LemA}

{\em Proof.}  First we will show that $H^n \subseteq \mbox{Ker}\psi_1$. 
Let $W = W(\mathcal A)$ be an $\mathcal A$-word. 
According to relations (5), (3), (6), we 
have in the group $G$ the  following equalities 
$$
1 \overset  G  = (W(\mathcal X) c W(\mathcal X)^{-1})^n \overset  G  =
( c W(\mathcal B))^n .
$$

Conjugating the last word by $y$ and using relations (4), (7), (9), we further have 
$$
1 \overset  G  = ( c W(\mathcal B))^n  \overset  G  =
 (y c W(\mathcal B) y^{-1})^n  \overset  G  = (c W(\mathcal B) W(\mathcal A) )^n .
$$
Making use of  relations (7)--(8), we can see that
$$
1 \overset  G  = ( c W(\mathcal B))^n W(\mathcal A)^n ,  
$$
whence $W(\mathcal A)^n  \overset  G    = 1$.  The inclusion  
$H^n \subseteq \mbox{Ker}\psi_1$ is proven.

To prove the converse, consider the following  presentation for 
$$
K = H/H^n =  \langle  a_i,    i \in I  \  \|\  \  R =1, \ R  \in \mathcal R, \ W^n =1, \ W \in F(\mathcal A)  \rangle 
$$
and let
$$
\gamma_0  : F(\mathcal A) \to  K
$$
 be the natural homomorphism. 

Let 
$$
B(\mathcal B\cup c, n) = F(\mathcal B\cup c)/  F(\mathcal B\cup c)^n
$$
be a free Burnside group of exponent $n$ in the alphabet $\mathcal B\cup c$
and 
$$
{\delta_0} : F(\mathcal B \cup c) \to B(\mathcal B \cup c, n)
$$
 be the natural homomorphism.  Consider the direct product
$$
P = B(\mathcal B \cup c, n) \times K .
$$
Note that the map
$$
c^{\delta_0} \to c^{\delta_0} b_i^{\delta_0},  \quad  
b_j^{\delta_0} \to b_j^{\delta_0} , \ j \in I, 
$$
extends to an automorphism of the Burnside group 
$B(\mathcal B \cup c, n)$.  Hence we can 
consider an  HNN-extension of $P$ with, say,  stable letter $x_i$ given by the following  relative presentation
$$
 \langle  \;  P,  \; x_i   \  \|\  \  
x_i c^{\delta_0} x_i^{-1} = c^{\delta_0} b_i^{\delta_0}  ,  \  
x_i b_j^{\delta_0} x_i^{-1} = b_j^{\delta_0} , \ j \in I \; \rangle  .
$$
Also note that the map 
$$
c^{\delta_0} \to c^{\delta_0}, \  \  
b_j^{\delta_0} \to b_j^{\delta_0} a_j^{\gamma_0}  , \  j \in I, 
$$
extends to an isomorphism of  $B(\mathcal B \cup c, n)$ to  a subgroup of $P$ generated by $c^{\delta_0}$,   $b_j^{\delta_0} a_j^{\gamma_0}$,  $j \in I$, because 
$K^n = \{ 1  \}$ and  $B(\mathcal B \cup c, n)$ is a free 
Burnside group of exponent $n$ with free generators 
$c^{\delta_0}$,   $b_j^{\delta_0}$, $j  \in I$.  

Therefore, we can consider a multiple  HNN-extension $E$ of $P$ with stable letters $y, x_i, i \in I,$ defined by the 
following  relative presentation
\begin{multline*}
E =  \langle \   P,  \; y,  \;  x_i, \ i \in I    \  \|\  \ 
x_i c^{\delta_0} x_i^{-1} = c^{\delta_0} b_i^{\delta_0} ,  \  
x_i b_j^{\delta_0} x_i^{-1} = b_j^{\delta_0} ,   \\ 
y c^{\delta_0} y^{-1} = c^{\delta_0} ,  \  
y b_k^{\delta_0} y^{-1} =  b_k^{\delta_0} a_k^{\gamma_0} ,  \  
\ i,j,k  \in I \   \rangle  .
\end{multline*}
It remains to notice that the natural homomorphism $F(\mathcal U) \to E$ factors out through $G$ and $K= F(\mathcal A)^{\gamma_0}$ naturally embeds in $E$. Therefore, the  converse inclusion
$\mbox{Ker}\psi_1 \subseteq H^n$  also holds and Lemma A is proven. \qed
\smallskip

Consider the quotient $G/G^{n^3}$ and let
$$
{\varepsilon_0} : F(\mathcal U) \to G/G^{n^3}
$$ 
denote the natural homomorphism.
 
\begin{LemB}  The map 
$a_i^{\beta_0} \to a_i^{\varepsilon_0}, i \in I$,  extends to a monomorphism
$$
\varphi_0 :  
F(\mathcal A)^{\beta_0}  \to  F(\mathcal U)^{\varepsilon_0} = G/G^{n^3} .
$$
\end{LemB}

{\em Proof.}  Clearly, the map $a_i^{\beta_0}  \to a_i^{\varepsilon_0}$, $i \in I$, 
extends to a homomorphism 
$$
\varphi_0 :  F(\mathcal A)^{\beta_0}  \to  F(\mathcal U)^{\varepsilon_0} 
$$
and we  have to show that Ker$\varphi_0 = \{ 1 \}$.
To do this we first repeat the construction of the proof  of Lemma A and define as there groups 
$$
K,  \ \   B(\mathcal B \cup c, n), \ \    P
$$ 
and natural  homomorphisms 
$$
{\gamma_0} :  F(\mathcal A)  \to K = H/H^n , \quad 
{\delta_0} :    F(\mathcal B \cup c)  \to B(\mathcal B \cup c, n)  .
$$
We also consider a free Burnside group  
$B(\mathcal X, n)$  in the alphabet 
$\mathcal X = \{ x_i \;  | \;  i \in I\}$ of exponent $n$ and let
$$
{\omega_0} : F(\mathcal X)  \to  B(\mathcal X, n)
$$
be the natural  homomorphism.  Note that  the map 
\begin{equation} \label{11}
c^{\delta_0} \to c^{\delta_0} b_i^{\delta_0}, \  b_j^{\delta_0} \to b_j^{\delta_0}, \  a_k^{\gamma_0} \to   a_k^{\gamma_0}, \ \  j, k \in I, 
\end{equation}
extends to an automorphism $\zeta_i, i \in I$, of the direct product 
$$
P = B( \mathcal B \cup c, n) \times K . 
$$ 
Moreover, it is easy to see that the subgroup 
$\langle \zeta_i \; | \;  i \in I \rangle \subseteq \mbox{Aut}P$ is isomorphic to a free Burnside group  of exponent $n$ which is 
freely generated by $\zeta_i, i \in I$.  
Hence, we can consider a semidirect product 
$$
Q = P \leftthreetimes  B(\mathcal X, n)
$$
of $P$ and  $B(\mathcal X, n)$ so that if $p \in P$ then
$$
x_i^{\omega_0} p x_i^{-{\omega_0}} = \zeta_i(p) .
$$
Clearly, the group $P$ has exponent $n$ and so the group $Q$ has exponent $n^2$.

Recall that, when proving Lemma A, we saw that the map
$$
{\varkappa_0} : c^{\delta_0} \to c^{\delta_0} , \   b_i^{\delta_0} \to b_i^{\delta_0} a_i^{\gamma_0},  \  i \in I, 
$$
extends to an isomorphism
\begin{equation} \label{12}
{\varkappa_0} :  \langle \; b_i^{\delta_0}, \; c^{\delta_0}  \  |  \ i \in I \; \rangle  \to   
\langle \; b_i^{\delta_0} a_i^{\gamma_0},  \;  c^{\delta_0} \  |   \ i \in I \;  \rangle
\end{equation}
of corresponding subgroups of $P$. Therefore, we can consider the following 
HNN-extension of  $Q$ with stable letter $y$ 
\begin{equation}  \label{13}
\mathcal G = \langle \; Q, y \ \|  \  
y c^{\delta_0} y^{-1} = c^{\delta_0} ,  \
y b_i^{\delta_0} y^{-1} = b_i^{\delta_0} a_i^{\gamma_0},  \   i \in I \;  \rangle . 
\end{equation}

The proof of  Lemma B will be continued in Sect. 3.

\section{Adjusting the Machinery of \cite{I94}}

Now we will make use of the machinery of  \cite{I94} to show (see Lemma C) 
that if $n > 2^{16}$ is odd then the quotient $\mathcal G/ \mathcal G^{n^3}$ 
of the group $\mathcal G$  defined by presentation (13) 
is infinite and $Q$ naturally embeds in  $\mathcal G/ \mathcal G^{n^3}$. Since the natural homomorphism
$$
K = H/H^n \to Q \to G /  G^{n^3}
$$
(naturally) factors out through the quotient $\mathcal G /  \mathcal G^{n^3}$, 
Lemma B will be proven. 

Diagrams over the group $\mathcal G(0) = \mathcal G$ given by 
presentation (13)  (or, briefly, 
over $\mathcal G$),   called  {\em diagrams of rank} 0, 
are defined to be maps that have two types of  2-cells.  

A 2-cell $\Pi$ of the first type, called a 0-{\em square}, has four edges in its counterclockwise oriented boundary
(called {\em contour})  
$$
\p\Pi = e_1 e_2 e_3 e_4
$$ 
and 
$$
\varphi(e_1) = \varphi(e_3)^{-1} = y , \ \   \varphi(e_2) = g,  \ \  \varphi(e_4) = {\varkappa_0}(g)^{-1} ,
$$
where, as in  \cite{I94},   $\varphi$ is the labeling function, 
$g \in \langle b_i^{\delta_0}, c^{\delta_0} \; | \;   i \in I \rangle \subseteq Q$
(perhaps, $g=1$) and $\varkappa_0$ is defined by (12). 
Observe that we use Greek letters with no indices exactly as in 
\cite{I94}  (in particular, see table (2.4) in \cite{I94}). 

A 2-cell $\Pi'$ of the second type, called a 0-{\em circle},  has $\ell \ge 2$ edges in its contour  $\p\Pi' = e_1 \dots e_\ell$ so that 
$\varphi(e_1), \dots ,  \varphi(e_\ell) \in Q$ and the word 
$$
\varphi(\p\Pi' ) = \varphi(e_1)  \dots  \varphi(e_\ell) 
$$ 
equals 1 in $Q$. 

Note that this definition of  2-cells in a diagram of rank 0 is a  generalization of corresponding definitions  in Ol'shanskii's book \cite{O89}.  According to the presentation (13), our basic alphabet is $Q_y = \{ Q,  y \}$ and, from now on 
(unless stated otherwise), all words will be those in the alphabet 
$Q_y^{\pm 1}= \{Q,  y, y^{-1}\}$, called $Q_y$-{\em words}.

Let 
$$
U_1  = S_{0,1} y^{\e_{1,1}} S_{1,1} \dots y^{\e_{\ell_1,1}} S_{\ell_1,1} , \quad
U_2  = S_{0,2} y^{\e_{1,2}} S_{1,2} \dots y^{\e_{\ell_2, 2}} S_{\ell_2,2} , 
$$
where $\e_{1,1}, \dots, \e_{\ell_1, 1},  \e_{1,2}, \dots, \e_{\ell_2, 2} \in \{\pm 1 \}$, 
$S_{0,1}, \dots,  S_{\ell_1,1}$ are $Q$-{\em syllables} of the word $U_1$ 
(that is, maximal subwords of $U_1$ all of whose letters are in $Q$; if 
$S_{0,1}$ is, in fact, missing in $U_1$,  then we assume that $S_{0,1}=1 \in Q$ and
set $S_{\ell_1,1}=1$ if $S_{\ell_1,1}$ is not present in $U_1$)  and 
$S_{0,2}, \dots, S_{\ell_2,2}$ are $Q$-syllables  of the word $U_2$. 

We will write $U_1 = U_2$ if $\ell_1 = \ell_2$, $\e_{j, 1}= \e_{j, 2}$  for all
$j =1, \dots, \ell_1$, and $S_{j, 1} = S_{j, 2}$ in the group $Q$ 
for all $j =0, \dots, \ell_1$. 

The {\em length}  $|U_1|$ of a word $U_1$ is $\ell_1$, that is, the number of occurrences of $y^{\pm 1}$ in $U_1$. 

If  the (images of)  words $U_1, U_2$ are equal in the group $\mathcal G(0) = \mathcal G$ given by (13)  then we will write $U_1 \overset 0  = U_2$.

By induction on $i  \ge 0$  we will construct groups $\mathcal G(i)$. Assume that the 
group $\mathcal G(i), i \ge 0$,  is already constructed as a quotient of the group of  $\mathcal G(0)$
by means of defining relations. Define $\mathcal X_{i+1}$ to be  a maximal set of all
$Q_y$-words of length $i+1$ (if any) with respect to the following three properties.

\begin{enumerate}
\item[(AB1)] Every word $A \in \mathcal X_{i+1}$ begins with $y$ or $y^{-1}$.
\item[(AB2)] The image of every word 
$A \in \mathcal X_{i+1}$ has infinite order in the group $\mathcal G(i)$. 
\item[(AB3)]  If $A, B$ are distinct elements of  $\mathcal X_{i+1}$  then the 
image of  $A^{n^2}$ is not conjugate in $\mathcal G(i)$ to the image of $B^{n^2}$ or 
$B^{-n^2}$. 
\end{enumerate}

Note  that it follows from the analogue of Lemma 18.2 in rank $i \geq 0$ that 
the set $\mathcal X_{i+1}$ is nonempty.

Similar to  \cite{O82}, \cite{O89}, \cite{I94}, we will call a word  $A \in \mathcal X_{i+1}$ a  {\em period of rank} $i+1$.

Now we define the group $\mathcal G(i+1)$  by imposing all relations 
$A^{n^3} =1$, called {\em  relations of rank $i+1$}, on the group  $\mathcal G(i)$:
$$
\mathcal G(i+1) = 
\langle \; \mathcal G(i) \; \| \;  A^{n^3} =1, \ A \in \mathcal X_{i+1} \;  \rangle .
$$

It is clear that 
$$
\mathcal G(i+1) = \langle \; \mathcal G  \; \| \;  A^{n^3} =1, \ A \in \bigcup_{j=1}^{i+1} \mathcal X_{j} \;  \rangle .
$$

We also define the limit group $\mathcal G(\infty)$  by imposing on the free product  
$Q * \langle y \rangle_\infty$ 
all relations of all ranks $i = 0, 1, 2, \dots$:
\begin{multline} \label{14}
\mathcal G(\infty) = \langle \  Q, \; y   \   \|   \  
y c^{\delta_0} y^{-1}=  c^{\delta_0}, \
y b_k^{\delta_0} y^{-1}=  b_k^{\delta_0} a_k^{\delta_0},  \    k \in I, \\
A^{n^3} =1, \  A \in \bigcup_{j=1}^{\infty} \mathcal X_{j} \  \rangle . 
\end{multline}

The main technical result relating to the group $\mathcal G(\infty)$ is the following 
(as was pointed out above, Lemma B is immediate from Lemma C).

\begin{LemC} Suppose that $n > 2^{16}$  is odd. Then the 
group $\mathcal G(\infty)$ given by the presentation $(14)$  has exponent $n^3$ 
and the group $Q$ naturally embeds in $\mathcal G(\infty)$.
\end{LemC}

{\em Proof.}
We will make use of the machinery of the article \cite{I94} to prove 
Lemma C (cf. \cite{I00}).  
In order to do that we will have to make necessary changes in definitions, statements of lemmas of \cite{I94}  and their proofs.  

First of all, we take into account that there are multiple periods of rank $i$ and the symbol $A_i$ will denote one of many periods of rank $i$ (note that the length $|A_i|$ of  $A_i$ is now $i$). Next, the number $n$ 
in \cite{I94}   is  replaced by $n^3$. 

In the definition of an $A$-periodic word, it is now assumed that $A$ starts with $y$ or $y^{-1}$ and $A$ is not conjugate in $\G$ 
to a power $B^\ell$ with $|B| < |A|$. 

In addition to cells of positive rank, 
we also have (as in \cite{O89}, \cite{IO96}) cells of rank 0 
(which are now either 0-squares or 0-circles). 
The equality $r(\D)= 0$ now means that all cells in $\D$ have rank 0.

If 
$$
p = e_1 \dots e_\ell , 
$$   where $e_1, \dots, e_\ell$ are edges,  
is a path in a diagram $\D$ of rank $i$ 
(that is, a diagram over the group  $\mathcal G(i)$) 
then the $y$-length $|p|$ of $p$ is $|\varphi(p)|$, that is, 
$|p|$ is the number of edges of $p$ labeled by $y^{\pm 1}$. 
The (strict) length of  $p$, that is,  the total number of edges of $p$,  
is $\ell$ and denoted by $\| p \|$.

In the definition (A1)--(A2)  of $j$-compatibility (p.13 \cite{I94}) we eliminate 
the part (A2)  because $n^3$ is odd and, 
similar to \cite{I94} in the case when $n$ is odd,     
it will be proven in a new version of Sect. 19 \cite{I94}  that 
there are no $\mathcal F(A_j)$-involutions and there is no  $j$-compatibility of 
type (A2).

We can also drop the definition of self-compatible cells (p.13 \cite{I94}) 
because they do not exist when $n^3$ is odd 
(which is again analogous to \cite{I94} in the case when $n$ is odd). Thus all lemmas in \cite{I94} whose conclusions deal with self-compatible cells, compatibility of  type (A2) actually claim that their assumptions are false (e.g.,  see lemmas of Sect. 12 \cite{I94}). On the other hand, the existence of  self-compatible cells in assumptions of lemmas of \cite{I94}  is now understood as the existence of  noncontractible  $y$-annuli 
which we are about to define. 

A $y$-{\em annulus} is defined to be an annular 
subdiagram $\G$ in a diagram $\D$ of rank $i$ 
such that $\G$ consists of 0-squares $S_1, \dots, S_k$ so that if
$$
\p S_\ell = f_{1, \ell} e_{1, \ell} f_{2, \ell} e_{2, \ell} ,
$$
where $e_{1, \ell}, e_{2, \ell}$ are $y$-edges 
(that is, labeled by $y^{\pm 1}$) of $\p S_\ell$,
$1 \le \ell \le k$, then  $e_{2, \ell}= e_{1, \ell+1}^{-1}$,  
where the second subscript is \ mod$k$, 
$\ell =1, \dots, k$. If  $\G$ is contractible into a point  in $\D$ then we will call
$\G$ a {\em contractible} $y$-annulus.  
Otherwise,  $\G$ is a  {\em noncontractible} $y$-annulus. 
If   $\G$ is contractible in $\D$ and bounds a simply connected subdiagram
$\G_0$ in $\D$ with $\varphi( \p \G_0) = 1$ in $Q$  then $\G$ is termed 
a  {\em reducible} $y$-annulus.

In the definition of a reduced (simply connected or not)  diagram $\D$ 
of rank $i$ (p.13 \cite{I94}),  we additionally require that $\D$ contain no  
reducible $y$-annuli. 

As in \cite{I94}, we can always remove reducible pairs and reducible $y$-annuli   
in a diagram $\D$ of rank $i$ to obtain from $\D$ a reduced  
diagram $\D'$ of rank $i$. Note that in general  it is not possible to get rid of noncontractible $y$-annuli  (in non-simply connected  diagrams of rank $i$). 

In the definition of a 0-bond $\E$  between $p$ and $q$ (p.15 \cite{I94}) we now require that
$\E$ consist of several  0-squares $S_1, \dots, S_k$ so that if 
$$
\p S_\ell = f_{1, \ell} e_{1, \ell} f_{2, \ell} e_{2, \ell} ,
$$
where $e_{1, \ell}, e_{2, \ell}$ are $y$-edges  of $\p S_\ell$,
$1 \le \ell \le k$, then $e_{1,1}^{-1} \in p$, $e_{2,\ell}= e_{1,\ell+1}^{-1}$, 
$\ell =1, \dots, k-1$, and $e_{2, k}^{-1} \in q$. 	

The {\em standard} contour of the 0-bond $\E$ between $p$ and $q$  is
$$
\p \E = e_{1, 1}^{-1}(f_{1, 1}^{-1} \dots f_{1, k}^{-1}) e_{2, k}^{-1}
(f_{2, k}^{-1} \dots f_{2, 1}^{-1})
$$
and the edges  
$e_{1,1}^{-1}, e_{2,k}^{-1}$ are denoted by $\E \wedge p$, $\E \wedge q$, 
respectively. 

In the definition of  a simple in rank $i$ word $A$ (p.19  \cite{I94}), 
we additionally require that $|A| >0$ and  $A$ start with $y^{\pm 1}$. 
Observe that it follows from Lemma 18.2 (in rank $i-1$)  and definitions that 
a period of rank $i$ is simple in rank $i-1$ (and hence in any rank $j \leq i-1$).

In the definition of a tame diagram of rank $i$ (p.19 \cite{I94}), we make two changes. First, in  property (D2), we require that if 0-squares $S_1, \dots, S_k$ form 
a subdiagram $\E$ as in  the definition of a 0-bond and $p= q = \p\Pi$, 
where $\Pi$ is a cell of rank $j >0$ in $\D$,  then 
$\E$ is a 0-bond between  $\p\Pi$ and $\p\Pi$ in $\D$. Second, we add the 
following property. 
\begin{enumerate}
\item[(D3)]  $\D$ contains no contractible $y$-annuli. 
\end{enumerate}

In the definition of  a complete system (p.23 \cite{I94}) we require in (E3) that $e$ be a 
$y$-edge. 

In Lemma 4.2, the strict length $\| s_1 \|,  \| s_2 \|$ of $s_1, s_2$ is meant. 

In the definition of  the weight function $\nu$ (p.28 \cite{I94}), 
we require in (F1) that $e$ be a $y$-edge.  
In (F2), we allow that $e$ is not  a $y$-edge.

In the beginning of the proof of Lemma 6.5, we note that the lemma is obvious if $\D$ contains no cells of positive rank. In general, repeating  arguments of \cite{I94}, we always understand "cells" as cells of positive rank and keep in mind the existence of cells of rank 0. 

In the definition of the height $h(W)$ of a word $W$ (p.89 \cite{I94}),  we additionally set 
$$
h(W) = \tfrac{1}{2}
$$  
if  $W \overset i  \neq 1$  and 
$W$ is conjugate in rank $i$ to a word $U$ with $|U| = 0$ 
(that is, $U \in Q$ and $U  \overset i  \neq 1$). 

In Lemma 10.2, we allow the extra case when $h(W) = \tfrac{1}{2}$. 

Here is a new version of Lemma 10.4.

\begin{Lem10.4}  $(a)$  If a word $W$ has finite order $d > 1$ in the group 
$\mathcal G(i)$ then $n^3$ is divisible by $d$. 

$(b)$ Every word $W$  with  $|W| \leq i$  has finite order in rank $i$.
\end{Lem10.4}

{\em Proof.} (a)  By  Lemma 10.2, either  $h(W) = \tfrac{1}{2}$ 
(and then our claim is immediate from the construction  of the group $Q$, 
see the definition of presentation (13)) or, otherwise, 
$W$ is conjugate in rank $i$ 
to a word of the form $A^k F$, where  $A$ is a period of rank $j \le i$,
$0 < k < n^3$  and $F \in \mathcal F(A)$.  
In the latter case, it follows from Lemma 18.5(c) in rank $j-1 <i$ that 
$$
(A^k F)^{n^2} \overset {j-1}  = A^{k n^2} .
$$
Therefore,
$$
(A^k F)^{n^3} \overset {j-1}  = A^{k n^3}  \overset {j}  = 1 ,
$$
whence $W^{n^3} \overset i = 1$ as required.  

(b)  By induction, it suffices to show that  every word $W$ with 
$|W| = i$ has finite order in rank $i$ (for $i=0$ this is obvious). 
It follows from the definition of periods 
of rank $i \geq 1$ that if $W$ has infinite order in rank $i-1$ then 
$W^{n^2}$ is conjugate in rank $i-1$ to $A^{\pm n^2}$, where 
$A$ is a period of rank $i$. 
Therefore, $W^{ n^3} \overset i = 1$ as desired. 

Lemma 10.4 is proved. \qed 
\smallskip

In Lemma 10.8, we drop part (b)  of its conclusion (and keep in mind that the term
"reducible cell" now means "$y$-annulus"). Note that the height
of $\varphi(q_1t_1)$ in Lemma 10.8 is $m \ge 1$ hence noncontractible $y$-annuli in $\D_0$, $\D_0^r$ are impossible (for otherwise, the height  $h(\varphi(q_1t_1))$ of  $\varphi(q_1t_1)$ would be at most $\tfrac {1}{2}$).

Lemma 10.9 is no longer needed 
for no path $q$ is (weakly) $j$-compatible with itself.

In the definition of a $U$-diagram of rank $i$ (p.134 \cite{I94}), we allow in property 
(U3) that the height $h(\varphi(e))$  of  $\varphi(e)$  is  $\tfrac {1}{2}$.

Lemma 12.3 now claims that there are no  $U$-diagrams of rank $i$. Recall that this agrees with our convention that if the conclusion  states the existence of self-compatible cells or $j$-compatibility of type (A2) then the assumption is false. 

The analogues of Lemmas 13.1--16.6 are not needed. 

In the hypothesis of  Lemma 17.1,  we now suppose that 
one can obtain  from $\D_0$ an annular reduced diagram 
of rank $i$  which contains no noncontractible $y$-annuli  by
means of  removal of  reducible pairs and reducible  $y$-annuli.  

According to our convention, in the statement of Lemma 17.2, we replace 
the  phrase "one has to remove a reducible cell  to reduce $\D_0$" by 
"one encounters a noncontractible $y$-annulus when reducing $\D_0$". 
In the conclusion of Lemma 17.2 and in its proof, we disregard reducible cells, 
$\mathcal F(A_j)$-involutions and consider, instead, noncontractible $y$-annuli and their 0-squares.  

The new version of Lemma 17.3 is stated as follows.

\begin{Lem17.3} Let $\D$ be a disk reduced diagram of rank $i$
whose contour is $\p\D =bpcq$, where $\varphi(p)$ and $\varphi(q)^{-1}$ are
$A$-periodic words and $A$ is a simple in rank $i$ word with $|A| = i+1$ 
(in particular, $A$ is a period of rank $i+1$).  
Suppose also that $\D$ itself is a contiguity
subdiagram between $p$ and $q$ with $\min(|p|,|q|) > L|A|$. 
Then, in the notation of Lemmas 9.1--9.4,   there
exists a rigid subdiagram of  the form $\D(m_1, m_2)$  in $\D$ $( = \D(1, k))$ 
such that
$$
r(\D(m_1, m_2))=0
$$
and the following analogues of  inequalities (17.25)--(17.26)  (p.222 \cite{I94}) hold.
\begin{align*}
|q(m_2,m_1)| & > |q(k,1)|-4.4|A| > (L-4.4)|A| , \\  
|p(m_1,m_2)| & > |p(1,k)|-4.4|A| > (L-4.4)|A|  , \\ 
|x_{m_1}|  & = |y_{m_2}| =0  . 
\end{align*}            
\end{Lem17.3}

{\em Proof.} To prove this new version of Lemma 17.3, 
we repeat the argument of the  beginning of the proof of Lemma 17.3 \cite{I94}. As there, 
making use of Lemma 17.2, we prove Lemma 17.3.1. After that, 
arguing as in the proof of  Lemma 17.3.2, it is easy to show, using Lemma 12.3, 
that $r(\D(m_1, m_2))=0$. \qed
\smallskip

When proving the analogue of  Lemma 18.2, we pick a word $B= B(a_1, a_2)$ 
in the  alphabet  $\{ a_1, a_2\}$ of length $i+1$ so that $B$ has the same 
properties as those  in \cite{I94} and, in addition, first and last letters of  
$B$ are distinct (the existence of such a word easily follows from Lemma 1.7 \cite{I94}).  
Next, consider a word
$B(a_1, a_2, q)$ in the alphabet $\{ a_1, a_2, Q\}$ which is obtained from
$B(a_1, a_2)$ by plugging in an element $q \in Q, q \neq 1$, between each pair 
of consecutive letters 
of the word $B(a_1, a_2)$. Then we replace each occurrence of the letter $a_1$ in 
$B(a_1, a_2, q)$ by $y$ and  each occurrence of the letter $a_2$ in 
$B(a_1, a_2, q)$ by $y^{-1}$. Clearly, we have a word $B(y, y^{-1}, q)$ with
$|B(y, y^{-1}, q)|= i+1$. Now, in view of Lemmas 10.2, 10.4, 
we can repeat the arguments of the proof of 
Lemma 18.2 without any changes. 
\smallskip

Let $A$ be a period of rank $i+1$.   By  
$$
\mathcal F(A)  
$$ 
denote a maximal subgroup of  $Q \subseteq \mathcal G(i)$
(the fact that $Q$ naturally embeds in $ \mathcal G(i)$ is immediate from 
Lemma 6.2)  with respect to the property that $A$ normalizes this subgroup
$\mathcal F(A) \subseteq \mathcal G(i)$.

Here is a new version of Lemma 18.3.

\begin{Lem18.3}  Suppose that $A$ is a period of rank 
$i+1$. Furthermore, let $\D$ be a disk reduced diagram of rank $i$  
such that  $\p\D = b p c q$, where
$\varphi(p), \; \varphi(q)^{-1}$ are $A$-periodic words with  
$\min(|p|,|q|) > 1.1 n^2 |A|$, and $\D$  itself  be a contiguity subdiagram 
between  sections $p$ and $q$.  
Then there exists a $0$-bond $\E$ in $\D$ with the standard contour 
$\p \E = b_\E p_\E c_\E q_\E$, where $(p_\E)_-$, $(q_\E)_+$ are  
phase vertices  of  $p, q$, respectively,  such that 
$$
\varphi(b_\E) \in  \mathcal F(A) ,  \quad 
\varphi(b_\E) \in  B(\mathcal B \cup c, n) \subset Q  ,  \quad 
 A^{n^2} \varphi(b_\E)  A^{-n^2} \overset i  =  \varphi(b_\E)  
$$ 
and for every  integer $k$ one has
$$
(A^k \varphi(b_\E))^{n^2}  \overset i  =  A^{k n^2 } .
$$ 
\end{Lem18.3}

{\em Proof.} Lemma 17.3 enables us to assume that 
$$
r(\D) = 0, \quad  \min( |p|, |q| ) > (1.1 n^2 - 5) |A| , \quad  |b| = |c|=0 .
$$
In particular, there are $| p|$ 0-bonds between $p$ and $q$ in $\D$. 
Let $\E$ be a 0-bond between sections $p$ and $q$ and 
$$
\p \E = b_\E p_\E c_\E q_\E 
$$
be the standard contour of $\E$, 
where $p_\E = \E \wedge p$, $q_\E = \E \wedge q$.  
It is clear that $\mbox{div}((p_\E)_-, (q_\E)_+)$  
does not depend on $\E$. Suppose that
\begin{equation} \label{15}
\mbox{div}((p_\E)_-, (q_\E)_+) \neq 0 .  
\end{equation} 

$A$-Periodically extending $p$ or $q^{-1}$ on the left as in the beginning of the  proof 
of Lemma 17.3 (see Fig. 17.4(a)--(b) in \cite{I94}), we will get a diagram $\D'$ with
$\p \D' = b' p' c' q'$ such that both $\varphi(p')$ and  $\varphi(q')^{-1}$ begin with a cyclic permutation $\bar A$ of $A$. 

As in the  proof of  Lemma 17.1,  we can easily get,  making use of (15),  that the annular diagram $\D'_0$ (obtained from $\D'$ as in Lemma 17.1) 
contains no $y$-annuli, in particular, $\D'_0$ is already reduced.  
Therefore,  Lemma 17.1 applies to $\D'$ and yields that
$\varphi(b') \overset i  = 1$. It follows from (15) that $|b'| >0$ and so $\bar A$ is not cyclically reduced in rank $i$. This contradiction proves that (15) is false and 
$$
\mbox{div}((p_\E)_-, (q_\E)_+)  =  0 . 
$$

Without loss of generality, we may assume that words 
$\varphi(p)$, $\varphi(q)^{-1}$ start with $A$ and  the word $A$ 
starts with $y^{-1}$.

Using the notation of Lemma 9.1, let $\E_1, \dots, E_{|p|}$ be all 
(consecutive along $p$) 0-bonds between $p$ and $q$ with standard contours
$$
\p \E_\ell = b_\ell  p_\ell  c_\ell  q_\ell , 
$$
where $p_\ell  = \E_\ell \wedge p$, $q_\ell  = \E_\ell \wedge q$ 
and  $1 \le \ell \le |p|$.

Also, we let $|A| = m$ and consider  words 
$$
V_0 = \varphi(b_1), \ V_1 = \varphi(b_{m+1}), \dots ,  
V_t =   \varphi(b_{mt+1}),  \dots ,  
V_{n^2 +1}=   \varphi(b_{m (n^2 +1) + 1}) 
$$
(recall that $|p| > (1.1n^2 - 5)|A| >  (n^2 + 5)|A|$).  
It is clear that
$$
A^{-1} V_t A = V_{t+1} 
$$ 
and $V_t, V_{t+1}$ are words in $Q$, 
where $t = 0, \dots ,  n^2$.

Recall that the group $Q = P  \leftthreetimes  B( \mathcal X \cup c, n)$  
is a semidirect product of the free Burnside group 
$$
B( \mathcal X, n) = \langle x_k^{\omega_0} \; |  \; k \in I \rangle
$$ 
and the direct product 
$$
P = K \times B( \mathcal B \cup c, n)
$$ 
of $K = H/H^n =  \langle a_k^{\gamma_0} \; | \;  k \in I \rangle$  and 
the free Burnside group  $B( \mathcal B \cup c, n) = \langle c^{\delta_0}, b_k^{\delta_0} \;  | \; k \in I \rangle$ 
so that the conjugation of $P$ by $x_k^{\omega_0}$ induces the automorphism (11).

Therefore, it follows from the definition of the group $\mathcal G = \mathcal G(0)$ 
given by (13) that 
every $V_t$, $t =0,1, \dots, n^2$, 
belongs to the subgroup  
$$
B( \mathcal B \cup c, n)  \subset Q
$$ 
which is  a factor of $P$  
(otherwise, $y V_t y^{-1} \not\in Q$, see (13)). 

Consider the quotient of $\mathcal G$ by the normal closure
$\langle  \langle H/H^n \rangle \rangle$ of 
$$
H/H^n =  \langle a_k^{\gamma_0} \; | \; k \in I  \rangle
$$
and let
$$
\varphi_0 : \mathcal G  \to   \mathcal G / \langle  \langle H/H^n \rangle \rangle
$$ 
be the natural homomorphism. 

It is clear that $Q^{\varphi_0}$ is naturally isomorphic to $Q$, 
so we can assume  that $Q^{\varphi_0}$ itself  
is a splitting extension  of
$B( \mathcal B \cup c, n) =  B( \mathcal B \cup c, n)^{\varphi_0}$ by 
$B( \mathcal X, n) = B( \mathcal X, n)^{\varphi_0}$ so that
$$
x_j^{\omega_0} c^{\delta_0} x_j^{-{\omega_0}} = c^{\delta_0} b_j^{\delta_0}, \ \  x_j^{\omega_0} b_{j'}^{\delta_0} x_j^{-{\omega_0}} =  b_{j'}^{\delta_0} , \ 
j, j' \in I .
$$

Observe that it follows from $V_t \in B( \mathcal B \cup c, n)$ that 
$V_t^{\varphi_0} = V_t$.  

Consider the quotient 
$\mathcal G^{\varphi_0} / \langle \langle  y^{\varphi_0} \rangle \rangle$ and let 
$$
{\varphi_1} : \mathcal G  \to \mathcal G^{\varphi_0} /
\langle \langle y^{\varphi_0} \rangle  \rangle
$$
denote the  natural homomorphism.  Clearly, 
$V_t^{\varphi_1} = V_t$, so 
$$
V_{t+1} = A^{-{\varphi_1}} V_t A^{\varphi_1} , \  \  t = 0,1, \dots, n^2-1 ,  
$$
and hence
$$
V_{n^2} = (A^{-{\varphi_1}})^{n^2} V_0 (A^{\varphi_1})^{n^2} . 
$$

Since $A^{\varphi_1} \in Q^{\varphi_1}$, 
$Q^{\varphi_0}$ and $Q^{\varphi_1}$ are naturally isomorphic  and 
$Q^{\varphi_0}$ has  exponent $n^2$, it follows  that
$V_{n^2} = V_0$. Now it is obvious that $V_0 \in \mathcal F(A)$ and
$$
A^{n^2} V_0 A^{-n^2}  \overset i  =  V_{n^2}  \overset i  =  V_0 .
$$
At last, we have 
\begin{multline*}
(A^k V_0)^{n^2}   \overset 0  =  
A^{k} V_0 A^{-k} A^{2k} V_0 A^{-2k} \dots 
A^{k n^2} V_0 A^{-k n^2}  A^{k n^2}  \overset i  = \\
(A^{\varphi_1})^{k} V_0 (A^{\varphi_1})^{-k} 
(A^{\varphi_1})^{2k} V_0 (A^{\varphi_1})^{-2k} \dots 
(A^{\varphi_1})^{k n^2} V_0 (A^{\varphi_1})^{-k n^2}  A^{k n^2}  \overset 0  = \\
((A^{\varphi_1})^k V_0)^{n^2} (A^{\varphi_1})^{-kn^2} A^{k n^2}  
 \overset i  = A^{k n^2} 
\end{multline*}
and Lemma 18.3 is proven. \qed
\smallskip

Lemma 18.4 is not needed. 

Let us state a new version of Lemma 18.5.

\begin{Lem18.5}
Let $A$ be a period of rank $i+1$ and  $\mathcal F(A)$ be a maximal subgroup of
$Q \subset \mathcal G(i)$ with respect to the property that $A$
normalizes $\mathcal F(A)$.  Then the following are true.

$(a)$  The subgroup $\mathcal F(A)$  is defined uniquely.

$(b)$   Suppose $\D$ is a disk reduced diagram of rank $i$ with $\p\D= bpcq$, 
where $\varphi(p)$, $\varphi(q)^{-1}$ are $A$-periodic words with
$\min(|p|,|q|) > \beta  n^3 |A|$, such that $\D$ itself is a contiguity subdiagram
between $p$ and $q$. Then there is a $0$-bond $\E$ between $p$ and $q$ 
with the standard contour $\p \E = b_\E p_\E c_\E  q_\E$, 
where $p_\E = \E \wedge p$, 
$q_\E = \E \wedge q$, such that $(p_\E)_-$,  $(q_\E)_+$ are phase vertices of $p, q$, respectively, and $\varphi(b_\E) \in \mathcal F(A)$.

$(c)$    $A^{n^2}$ centralizes  the subgroup $\mathcal F(A)$ and 
$\mathcal F(A)$ has  exponent $n$.   Furthermore, 
if $F \in \mathcal F(A)$ and $k$ is an  integer  
then 
$$
(A^k F)^{n^2} \overset i  = A^{k n^2} .
$$

$(d)$   The subgroup $\langle \mathcal F(A), \; A \rangle$ of $\mathcal G(i)$ has the property that
a word $X$ belongs to $\langle \mathcal F(A), \; A \rangle$
if and only if there is an integer $m \ne 0$ such that 
$X A^m X^{-1} \overset i = A^m$. 
\end{Lem18.5}

{\em Proof.} (a)  This is obvious from definitions.

(b)  Since $\beta  n^3 |A| > 1.1n^2 |A|$,  we can apply Lemma 18.3 
to $\D(m_1, m_2)$  which yields the required conclusion. 

(c) Suppose that $X \in \mathcal F(A)$. Picking sufficiently large $\ell$, we have 
an equality
$$
A^{-\ell} X A^\ell  \overset i  = Y ,
$$
where $Y \in  \mathcal F(A)$.  Consider a reduced diagram $\D$ 
of rank $i$ for this equality. 
Using Lemmas 6.5, 3.1 and inequality $\ell \gg 1$, 
it is easy to show  that $\D$ contains a contiguity subdiagram
$\D_0$ between sections whose labels are $A^{-\ell}$ and $A^\ell $ with the standard
contour
$$
\p \D_0 = b p c q ,
$$ 
where $\varphi(p)$,  $\varphi(q)^{-1}$ are $A$-periodic sections longer than 
${\beta} n^3 |A|$. 
Now we can refer to Lemma 18.3 and conclude that $X$ belongs to the double coset
$$
\langle A \rangle \varphi(b_\E)  \langle A \rangle \subseteq \mathcal G(i) , 
$$ 
where $\varphi(b_\E)$ is from the conclusion of Lemma 18.3. 
Since $X \in Q$ has finite order, it follows from Lemma 18.3 that 
$X  \overset i  =  A^{\ell'}  \varphi(b_\E) A^{-\ell'}$ with some $\ell'$. 
Therefore, $X$ commutes with $A^{n^2}$ in rank 
$i$  just like  $\varphi(b_\E)$  does and, by Lemma 18.3,   
$$
X^n  \overset  i    = 1, \quad 
(A^{k} X)^{n^2} \overset i  = A^{kn^2} .
$$ 
Part (c) is proven.

(d) By part (c), it suffices to show that an equality
$$
X  A^{m} X^{-1} \overset i  = A^{m} , 
$$
where $m \neq 0$,  implies that $X \in  \langle A,   \mathcal F(A) \rangle$. 
Arguing  exactly as in the proof of part (c), we can show that 
$$
X \in \langle A \rangle \varphi(b_\E)  \langle A \rangle \subseteq \mathcal G(i) , 
$$ 
where $\varphi(b_\E) \in  \mathcal F(A)$ by Lemma 18.3. Thus 
$X \in  \mathcal F(A)$  and Lemma 18.5 is proven.  \qed 
\smallskip

Here is a new version of Lemma 19.1.

\begin{Lem19.1}  There is no disk diagram  $\D$ of rank $i$ such that 
$\p \D = b p c q$, where $p$, $q$ are
$A$-periodic sections with $|p|, |q|  > \tfrac 13 {\beta} n^3 |A|$, $A$ is a period of rank $i+1$, and $\D$ itself is a contiguity subdiagram between $p$ and $q$. 
\end{Lem19.1}

{\em Proof.} Arguing on the contrary, we assume the existence of such a diagram $\D$
and,  replacing the coefficient $N$  ($N$ is defined in  (17.1) on p.212  \cite{I94}) 
by $\tfrac 13 {\beta} n^3$, repeat the proof of 
Lemma 19.1 \cite{I94}  up to getting equality (19.23) (p.290 \cite{I94}) which now reads
\begin{equation} \label{16}
\varphi(d) A^{n^2} \varphi(d)^{-1} \overset i  =  A^{-n^2} . 
\end{equation} 

\begin{Lem19.1.3}  The equality $(16)$ is impossible. 
\end{Lem19.1.3}

{\em Proof.} Arguing on the contrary, we note that it follows from (16) that 
\begin{equation} \label{17}
\varphi(d)^2 A^{n^2} \varphi(d)^{-2} \overset i  =  A^{n^2} .  
\end{equation}  

Recall that $|A| = i+1$ and, by Lemma 19.1.1, $\varphi(d)$ is conjugate  in rank $i$
to a word $W$ with 
$$
|W| < (0.5+\xi )|A| <|A| . 
$$ 
Hence, by Lemma 10.4(b),   $\varphi(d)$ has finite order  in rank $i$ and, 
by Lemma 10.4(a),   
$$
\varphi(d)^{n^3} \overset i = 1 .  
$$
This means  that 
$\varphi(d)  \in \langle \varphi(d)^2 \rangle$ in rank $i$ and so equalities (16)--(17) imply that 
$$
A^{2n^2}   \overset i  = 1 .
$$ 

A contradiction  to the definition of a period of rank $i+1$ completes  proofs 
of  Lemmas 19.1.3 and 19.1. \qed
\smallskip

Analogues of Lemmas 19.2--19.6 are no longer needed.

The statements and proofs of  Lemmas 20.1--20.2 are retained.
\smallskip 

Having made all necessary changes, we can now turn to the group $\mathcal G(\infty)$ given  by the presentation (14). 

It follows from Lemma 10.4(b) 
that every word $W$ has finite order in rank $i  \geq |W|$. Then, by 
Lemma 10.4(a), 
$$W^{n^3}  \overset i  = 1$$ 
provided that $i \geq |W|$. 
Thus, the group $\mathcal G(\infty)$ has exponent $n^3$. 

Suppose that $W$ is a word with $|W| = 0$, that is, $W \in Q$. Let $W = 1$ in 
the group $\mathcal G(\infty)$. Then there is an $i$ such that $W \overset i  = 1$.
Consider  a reduced diagram $\D$ for this equality. Since $|\p \D | =0$, it follows from 
Lemma 6.2 that $r(\D) =0$, that is, $\varphi(\p \D) =1$ in $Q$. 
Thus $Q$ naturally embeds in $\mathcal G(\infty)$.

Lemmas B and C are proven. \qed

\section{Proofs of Theorems A and B}

We start with proving Theorem B. Let $N_0  > 2^{16}$ be odd and 
$$
G_0 = \langle \; \mathcal A_0 \;  \|  \; R=1, R \in \mathcal R_0 \; \rangle
$$ 
be a finite presentation of a group $G_0$ with $|  \mathcal A_0 | > 1$.

Applying the construction of  presentation (10), we set $\mathcal A = \mathcal A_0$, 
$\mathcal R = \mathcal R_0$,  $H = G_0$, $n = N_0$. 
Denote the corresponding sets of
generators and relators  of  thus constructed presentation (10)  by 
$\mathcal A_1, \mathcal R_1$.  The group given by (10) denote by $G_1$. 
Since $\mathcal A_0 \subset \mathcal A_1$, $\mathcal R_0 \subset \mathcal R_1$, 
we can speak of  the  subgroup $\langle \mathcal A_0 \rangle$ of $G_1$ 
generated by the set $\mathcal A_0$.  
It follows from Lemma A that  the  subgroup 
$\langle \mathcal A_0 \rangle \subseteq  G_1$   
is naturally isomorphic to  the quotient $G_0/G_0^{N_0}$. 

At the second step, we again apply the construction of (10) with
$\mathcal A = \mathcal A_1$,  $\mathcal R = \mathcal R_1$,  $H = G_1$, and $n = N_1 = N_0^3$.
Denote the corresponding sets of generators and relators  of  
thus constructed presentation (10) by  $\mathcal A_2$, $\mathcal R_2$ and denote the group given 
by (10) by $G_2$.  
It follows from Lemmas A--B that the subgroups $\langle \mathcal A_0 \rangle$,  
$\langle \mathcal A_1\rangle$
of $G_2$  are  naturally isomorphic to 
$G_0/G_0^{N_0}$,  $G_1/G_1^{N_1}$,  respectively, and 
$G_0 / G_0^{N_0}$  naturally embeds in $G_1/G_1^{N_1}$. 

Keeping on doing this, we will inductively define a finite presentation 
$$
G_{\ell+1} = \langle \; \mathcal A_{\ell+1} \;  \|  \; R=1, \  
R \in \mathcal R_{\ell+1} \; \rangle ,
\quad \ell \geq 1 ,
$$ 
which is constructed as  presentation (10), where 
$\mathcal A = \mathcal A_\ell$,  $\mathcal R = \mathcal R_\ell$,  
$H = G_\ell$, and 
$$
n = N_\ell =  N_{\ell -1}^3=   N_0^{3^{\ell}} .
$$
It  follows from Lemmas A--B that  the subgroups $\langle \mathcal A_{\ell-1}\rangle$, 
$\langle \mathcal A_\ell \rangle$ of  $G_{\ell+1}$  are  naturally isomorphic to 
$$
G_{\ell-1}/G_{\ell-1}^{N_{\ell-1}} \subseteq 
G_\ell/G_\ell^{N_\ell}, \qquad G_\ell/G_\ell^{N_\ell} ,
$$
respectively. Therefore, by induction, the subgroup 
$\langle \mathcal A_k \rangle$ of $G_{\ell+1}$, where $0 \le k \le \ell$, is  naturally isomorphic to 
$G_k /G_k^{N_0^{3^k}}$.

Now define
$$
G_{\infty} = \langle \ \mathcal A_{\infty}= \bigcup_{j=0}^\infty \mathcal A_j 
 \  \|  \ R=1, \  R \in  \mathcal R_{\infty} = \bigcup_{j=0}^\infty \mathcal R_j   \ \rangle .
$$

It is clear that for every $\ell \ge 0$ the subgroup $\langle \mathcal A_\ell \rangle$ of 
$G_{\infty}$ is  naturally isomorphic to  $G_\ell/ G_\ell^{N_\ell}$  and that 
if $W_1, \dots, W_m$ are words in 
$ \mathcal A_\ell^{\pm 1}$ then they generate 
naturally isomorphic 
subgroups in $G_\infty$ and in $G_{\ell+1}$. Note that it is immediate from
the construction of the group $\mathcal G(\infty)$ 
(see (14) in which we put  
$\mathcal A = \mathcal A_{\ell-1}$,  
$\mathcal R = \mathcal R_{\ell-1}$,  $H = G_{\ell-1}$,
$n =  N_0^{3^{\ell -1}}$)  and  Lemma C that
$G_\ell/ G_\ell^{N_0^{3^\ell}}$ is not a locally finite group when $\ell \geq 1$. 
The proof of Theorem B is complete. 
\smallskip

Theorem A is immediate from Theorem B if we pick a free group of rank 2 as $G_0$
(then  $G_0/G_0^{N_0}$ is a free 2-generator  Burnside  group $B(2, N_0)$
of exponent $N_0$).

\end{document}